\documentclass[reqno]{amsart}
\usepackage{amsfonts,amsmath,amssymb}

\numberwithin{equation}{section}

\newtheorem{theorem}{Theorem}[section]
\newtheorem{lemma}{Lemma}[section]

\newtheorem{corollary}{Corollary}[section]
\newtheorem{proposition}{Proposition}[section]

\begin{document}
\title[Scattering for Beam in small dimensions]
{Scattering for the Beam equation in low dimensions}
\author{Benoit Pausader}
\address{Mathematics Department,
Box 1917,
Brown University,
Providence, RI 02912}
\email{Benoit.Pausader@math.brown.edu
%\\wstrauss@math.brown.edu
}

\date{January 16, 2009}

\begin{abstract}

In this paper, we prove scattering for the defocusing Beam equation $u_{tt}+\Delta^2u+mu+\lambda\vert u\vert^{p-1}u=0$ in low dimensions $2\le n\le 4$ for $p>1+8/n$. The main difficulty is the absence of a Morawetz-type estimate and of a Galilean transformation in order to be able to control the Momentum vector. We overcome the former by using a strategy of Kenig and Merle \cite{KenMer} derived from concentration-compactness ideas, and the latter by considering a Virial-type identity in the direction orthogonal to the Momentum vector.
\end{abstract}

\maketitle

The goal of the present paper is to further develop the scattering theory for the nonlinear beam equation
\begin{equation}\label{Beam}
\partial_t^2u+\Delta^2u+mu+\lambda\vert u\vert^{p-1}u=0,
\end{equation}
where $u:\mathbb{R}\times\mathbb{R}^n\to\mathbb{R}$, $1\le n\le 4$, $m>0$ and $1+8/n<p<+\infty$. Scattering for the Beam equation \eqref{Beam} was conjectured in Levandosky and Strauss \cite{LevStr} when the authors found an analogue of the Morawetz estimate for the Klein-Gordon equation that applied to the Beam equation. We refer also to the work of
Levandosky \cite{Lev1,Lev2} where scattering for small initial data was proved and an analysis of existence and stability
of traveling waves for (0.1) was carried over. We also mention Hebey and Pausader \cite{HebPau} for other results in the focusing case.

In Pausader \cite{PauBeam} the author proved scattering for the Beam equation in the high dimensional setting $n\ge 5$ where one can use the Morawetz estimate of Levandosky and Strauss \cite{LevStr}, the fast decay of the linear propagator and ideas from Tao \cite{Tao}. We also refer to Pausader and Strauss \cite{PauStr} for a discussion of various properties of the scattering operator.

Lower dimensions are more difficult to handle because linear solutions do not decay sufficiently fast and no analogue of the Morawetz estimate is known. Therefore, the method of Pausader \cite{PauBeam} completely breaks down.
In this paper, we use a different method to prove another part of the conjecture, namely we prove that scattering also holds for the low dimensions $2\le n\le 4$.
We also prove that, in dimension $n\ge 2$, any nonlinear solution has bounded space-time norms, with a bound depending only on the energy of the initial data.

\medskip

Our main result is as follows. See the next section for the definition of strong solution and of the energy.

\begin{theorem}\label{MainThm}
Let $2\le n\le 4$, $\lambda>0$ and $1+8/n<p<+\infty$. Then any strong solution $u$ of \eqref{Beam} scatters in the sense that there exists two linear solutions of the associated linear problem $v^\pm$ such that
$$\Vert u(t)-v^\pm(t)\Vert_{H^2}+\Vert u_t(t)-v^\pm_t(t)\Vert_{L^2}\to 0$$
as $t\to\pm\infty$.
\end{theorem}
It is easy to see (cf Proposition \ref{PropositionVersion6} below) that, when $1+8/n<p<+\infty$, and $n\le 4$, scattering for the Beam equation is implied by the finiteness of the following norm
\begin{equation}\label{DefOfNorm}
\Vert u\Vert_{S(I)}=\Vert u\Vert_{L^{p+1}(I\times\mathbb{R}^n)}
\end{equation}
when $I=\mathbb{R}$.
We actually prove Theorem \ref{MainThm} by proving, see Theorem \ref{Thm2} in Section \ref{CritSol},
that the $S$-norm is not only finite for all strong solutions, but also bounded by a quantity depending only on the energy.

\medskip

In order to prove Theorem \ref{MainThm}, we use a strategy derived from concentration-compactness ideas devised by Kenig and Merle \cite{KenMer,KenMer2} in order to tackle scattering for focusing critical problems. We also refer to Duyckaerts, Holmer and Roudenko \cite{DuyHolRou} for another adaptation in the subcritical case. As an important remark, it should
be noted that the Beam equation \eqref{Beam} does not seem to possess an invariance that could play the role of the Galilean invariance in setting the momentum of the critical solution to $0$, contrary to the Schr\"odinger equation.
We also mention that our approach is
somewhat different from the approach developed in the second order case by
Nakanishi
 \cite{Nak} and Colliander, Grillakis and Tzirakis \cite{ColGriTzi} where the authors prove scattering for Schr\"odinger and wave equations in small dimensions using an adaptation of the Morawetz inequality which is not available in the case of the Beam equation. Indeed, the lack of these Morawetz-type inequalities and of a Galilean-type invariance are major impediments in the case of the Beam equation.

More precisely, we proceed as follow. We define the function $\Lambda$ in Theorem \ref{Thm2} by
\begin{equation}\label{DefOfLambda}
\Lambda(E)=\sup\{\Vert u\Vert_{S(\mathbb{R})}^{p+1}:E(u)\le E\}
\end{equation}
where the supremum is taken over all nonlinear solutions of \eqref{Beam} of energy less than $E$, and
\begin{equation}\label{DefOfEmax}
E_{max}=\sup\{E:\Lambda(E)<+\infty\}.
\end{equation}
Our goal in the following is to prove that $E_{max}=+\infty$.
Using a profile decomposition similar to the one of Bahouri and G\'erard \cite{BahGer}, and a strategy inspired by Kenig and Merle \cite{KenMer} and Tao, Visan and Zhang \cite{TaoVisZha}, we prove here that if $E_{max}<\infty$, then there exists a nonlinear solutions of \eqref{Beam} of energy exactly $E_{max}$ with strong compactness properties. Then, using a Virial-type estimate, we find a contradiction. The Morawetz/Virial functional we use in the last step is inspired by a proof of nonexistence of traveling waves for the defocusing equation. This proof breaks down in dimension $n=1$, when a directional derivative is equal to a full derivative. This is the reason why we cannot conclude in the case $n=1$. We refer to Section \ref{ProofofThm} for more details.

This paper is organised as follows: in Section \ref{Not} we give notations we use throughout the paper, as well as previous results from Pausader \cite{PauBeam}. In Section \ref{CritSol}, we establish a profile decomposition for linear solutions associated to \eqref{Beam} and use it to produce a critical solution which is essentially a traveling wave, and in Section \ref{ProofofThm}, we finish the proof of Theorems \ref{MainThm} and \ref{Thm2}. Finally, in the appendix, we give two results about the translation parameter of the critical equation.

\section{Notations and previous results}\label{Not}

We fix notations we use throughout the paper. In what follows, we write $A\lesssim B$ to signify that there exists a constant $C$ depending only on $n$ such that $A\le C B$. When the constant $C$ depends on other parameters, we indicate this by a subscript, for exemple, $A\lesssim_u B$ means that the constant may depend on $u$. Similar notations hold for $\gtrsim$.  Similarly we write $A\simeq B$ when $A\lesssim B\lesssim A$.

\medskip

We assume that $n\le 4$.
In the sequel, we let $\mathcal{E}=H^2\times L^2$ be the energy space. On $\mathcal{E}$, we define the inner product
$$\langle (u_0,u_1),(v_0,v_1)\rangle_{\mathcal{E}}=\int_{\mathbb{R}^n}\left(u_1v_1+\Delta u_0\Delta v_0+ mu_0v_0\right)dx,$$
the free energy
$$E_0(u,v)=\frac{1}{2}\int_{\mathbb{R}^n}\left(v^2+(\Delta u)^2+mu^2\right)dx,$$ and we choose $\sqrt{E_0}$ as norm on $\mathcal{E}$.
We also define two other functionals on $\mathcal{E}$, the nonlinear energy,
\begin{equation}\label{SecNotEnengy}
E(u,v)=E_0(u,v)+\frac{\lambda}{p+1}\Vert u\Vert_{L^{p+1}}^{p+1}
\end{equation}
and the momemtum,
\begin{equation}\label{SecNotMom}
\hbox{Mom}(u,v)=\int_{\mathbb{R}^n}v\nabla udx.
\end{equation}
Given $(u_0,u_1)\in\mathcal{E}$, there exists a unique solution $w\in C(\mathbb{R},\mathcal{E})$ of the linear equation
\begin{equation}\label{LinearBeam}
w_{tt}+\Delta^2w+mw=0
\end{equation}
such that $(w(0),w_t(0))=(u_0,u_1)$.
We define the free evolution isometry group on $\mathcal{E}$, $\mathcal{W}$, by  $(w(t),w_t(t))=\mathcal{W}(t)(u_0,u_1)$. We also let $\pi_1$ and $\pi_2$ be the projection from $\mathcal{E}$ onto the first and second argument.
For $w$ a linear solution of \eqref{LinearBeam}, the free energy is constant, and we may write $E_0(w)$ instead of $E_0(w(t),w_t(t))=E(w(0),w_t(0))$, where $t\in\mathbb{R}$ is any time.

For $y\in\mathbb{R}^n$ and $f$ a function defined on $\mathbb{R}^n$, we define $\tau_yf$ by
\begin{equation}\label{SecNotTau}
\tau_yf(x)=f(x-y).
\end{equation}
For $1\le p,q<\infty$, we define $L^p(\mathbb{R},L^q)$ to be the completion of $C^\infty_c(\mathbb{R}^n)$ with the norm
$$\Vert f\Vert_{L^p(L^q)}=\left(\int_{\mathbb{R}}\left(\int_{\mathbb{R}^n}\vert f(t,x)\vert^qdx\right)^\frac{p}{q}dt\right)^\frac{1}{p},$$
and for $I$ an interval, we let $L^p(I,L^q)$ be the restiction to $I$ of functions in $L^p(\mathbb{R},L^q)$, and we endow $L^p(I,L^q)$ with the natural restriction norm. We adopt similar definitions when $p$ or $q$ is infinite. Given
$a\geq 1$, we let $a^\prime$ be the conjugate of $a$, so that
$\frac{1}{a}+\frac{1}{a^\prime}=1$.

We also let
$\mathbb{E}_I=L^\infty(I,H^2)\cap W^{1,\infty}(I,L^2)\cap W^{2,\infty}(I,H^{-2})\cap L^{p+1}(I,L^{p+1})$. We call strong solution of \eqref{Beam} on $I$ any function $u$ such that for all compact $J\subset I$, $u\in\mathbb{E}_J$ and \eqref{Beam} holds in $H^{-2}$ on $J$. It follows from results in Pausader \cite{PauBeam} (see also Levandosky \cite{Lev1,Lev2} for earlier results) that for any initial data $(u_0,u_1)\in \mathcal{E}$, there exists a unique nonlinear strong solution of \eqref{Beam}, $u\in C(\mathbb{R},H^2)\cap C^1(\mathbb{R},L^2)$ with initial data $(u(0),u_t(0))=(u_0,u_1)$. Besides, that solution has conserved energy and momentum, that we simply refer to as $E(u)$ and $\hbox{Mom}(u)$, thus for all $t$,
\begin{equation*}
\begin{split}
&E(u)=E(u(t),u_t(t))=E(u_0,u_1),\hskip.1cm\hbox{and}\\
&\hbox{Mom}(u)=\hbox{Mom}(u(t),u_t(t))=\hbox{Mom}(u_0,u_1),
\end{split}
\end{equation*}
where we define the energy and momentum on $\mathcal{E}$ by \eqref{SecNotEnengy} and \eqref{SecNotMom}.

For convenience, given a function $\omega\in\mathbb{E}_I$, we may write $E(\omega(t))$ or $E_0(\omega(t))$ instead of $E(\omega(t),\omega_t(t))$ and $E_0(\omega(t),\omega_t(t))$. We may also use the vectorial notation $\bar{\omega}(t)=(\omega(t),\partial_t\omega(t))$. In addition to the $S$ norm defined in \eqref{DefOfNorm},
we also need the following norm
$$\Vert u\Vert_{N(I)}=\Vert u\Vert_{L^\frac{2(n+2)}{n+4}(I,L^\frac{2(n+2)}{n+4})}+\Vert u\Vert_{L^\frac{2(n+4)}{n+8}(I,L^\frac{2(n+4)}{n+8})}.$$
We may omit $I$ when $I=\mathbb{R}$.

The Strichartz estimates from Pausader \cite{PauBeam} implies the following estimate when $n\le 4$:
let $u\in C(I,H^{-2})$ satisfy the equation $u_{tt}+\Delta^2u+mu=h$ with initial data $(u(0),u_t(0))\in\mathcal{E}$, and let $(a,b)$ be such that
\begin{equation*}
\begin{split}
2<a,b<\infty\hskip.1cm\hbox{and}\hskip.1cm \frac{4}{a}+\frac{n}{b}\le \frac{n}{2},
\end{split}
\end{equation*}
then
\begin{equation}\label{Stric}
\Vert u\Vert_{\mathbb{E}_I}+\Vert u\Vert_{L^a(I,L^b)}\lesssim \Vert u_0\Vert_{H^2}+\Vert u_1\Vert_{L^2}+\Vert h\Vert_{N(I)}.
\end{equation}
We note also that we can replace the control on $h$ in the $N$-norm by any term like
\begin{equation}\label{SecNotAlternativeNorms}
\Vert h\Vert_{L^{q^\prime}(L^{r^\prime})}+\Vert h\Vert_{L^{s^\prime}(L^{\sigma^\prime})},
\end{equation}
where $2\le q,r,s,\sigma< \infty$ and $4/q+n/r=2/s+n/\sigma=n/2$.

Finally, we need the following Littlewood-Paley projection.
Let $\psi\in C^\infty_c(\mathbb{R}^n)$ be supported in the ball
$B(0,2)$, and such that $\psi=1$ in $B(0,1)$. For any dyadic
number $N=2^k,k\in\mathbb{Z}$, we define the following
Littlewood-Paley operators:
\begin{equation}\label{DefLitPalOp}
\begin{split}
&\widehat{P_{\leq N}f}(\xi)=\psi(\xi/N)\hat{f}(\xi),\\
&\widehat{P_{>N}f}(\xi)=(1-\psi(\xi/N))\hat{f}(\xi),\\
&\widehat{P_Nf}(\xi)=\left(\psi(\xi/N)-\psi(2\xi/N)\right)\hat{f}(\xi).
\end{split}
\end{equation}
Similarly we define $P_{<N}$ and $P_{\ge N}$ by the equations
$$P_{<N} = P_{\leq N}-P_N\hskip.2cm\hbox{and}\hskip.2cm P_{\ge N} = P_{> N} + P_N,$$
and we let these operators act on couples of function as follows:
$P_N(u,v)=(P_Nu,P_Nv)$, and similarly for $P_{>N}$ and $P_{<N}$.
These operators commute one with another. They also commute with
derivative operators and with the semigroup
$\mathcal{W}$. In addition they are
self-adjoint and bounded on $L^q$ for all $1\le q\le\infty$.
Moreover, they enjoy the following Bernstein property:
\begin{equation}\label{BernSobProp}
\begin{split}
&\hskip.8cm \Vert P_{\ge N}f\Vert_{L^q}\lesssim_s
N^{-s}\Vert\vert\nabla\vert^sP_{\geq N}f\Vert_{L^q}\lesssim_s N^{-s}\Vert\vert\nabla\vert^sf\Vert_{L^q}\\
&\hskip.8cm\Vert\vert\nabla\vert^sP_{\le N}f\Vert_{L^q}\lesssim_s N^s\Vert P_{\le
N} f\Vert_{L^q}
\lesssim_s N^s\Vert f\Vert_{L^q}\\
&\hskip.8cm\Vert \vert\nabla\vert^{\pm s}P_Nf\Vert_{L^q}\lesssim_s
N^{\pm s}\Vert P_Nf\Vert_{L^q}\lesssim_s N^{\pm s}\Vert f\Vert_{L^q}\\
\end{split}
\end{equation}
for all $s\ge 0$, and all $1\le q\le\infty$, independently of $f$, $N$, and $q$, where
$\vert\nabla\vert^s$ is the classical fractional differentiation
operator. We refer to Tao \cite{TaoBook} for more details. 

\medskip

Now, we state previous results about solutions of \eqref{Beam}. In the following, we always assume that $1\le n\le 4$ and that $1+8/n<p<\infty$ or that $n\ge 5$ and $1+8/n<p<1+8/(n-4)$, and in this case, we replace $\Vert u\Vert_S$ by $\Vert u\Vert_{S^*}$.
A first result we need is the following lemma that allows us to understand when the linear solution is a good approximation of the nonlinear one.

\begin{lemma}\label{CondForScat}
Let $(u_0,u_1)\in\mathcal{E}$ and let $(w,w_t)=\mathcal{W}(\cdot)(u_0,u_1)$. There exists $\delta>0$ depending only on $\vert\lambda\vert$, $p$, $n$ and $E$ such that
if $E(u_0,u_1)\le E$ and $I$ is an interval such that
\begin{equation}\label{LWPStatement}
\Vert w\Vert_{S(I)}\le \delta
\end{equation}
then, letting $u$ be the strong solution of \eqref{Beam} with initial data $(u_0,u_1)$, we have that $u$ exists on $I$ and that
\begin{equation}\label{ConsLWP}
\Vert u\Vert_{S(I)}\lesssim \Vert w\Vert_{S(I)}.
\end{equation}
Besides, if $I=[T,+\infty)$, then $u$ scatters at $+\infty$.
In particular, $\Lambda$ in \eqref{DefOfLambda} is finite and sublinear in a neighborhood of $0$.
\end{lemma}
The proof follows by a straighforward fixed-point argument.
This proves that the Cauchy problem is well-posed in $C([0,T],H^2)\cap C^1([0,T],L^2)\cap L^{p+1}([0,T],L^{p+1})$ for initial data $(u_0,u_1)\in\mathcal{E}$ belonging to a ball in the energy space, where $T$ depends on the radius of the ball.
Using the conservation of energy, these solutions can be extended globally. This gives the following proposition which is the starting point of our investigation.

\begin{proposition}\label{SecNotGlobalExistenceLocallyBoundedProp}

For all initial data $(u_0,u_1)$ in the energy space, there exists a unique globally defined nonlinear solution $u\in C(\mathbb{R},H^2)\cap C^1(\mathbb{R},L^2)$. Besides the evolution flow $(u(0),u_t(0))\in\mathcal{E}\mapsto u\in \mathbb{E}_I$ is continuous for all compact time interval $I$.
\end{proposition}

The key point to understand better the behaviour of these solutions is to gain access to global in time bounds. Actually, as previously said, only a bound on the $S$ norm is sufficient. This is the object of the following proposition.

\begin{proposition}\label{PropositionVersion6}
Assume that $\lambda>0$, $1\le n\le 4$ and $p\ge (n+8)/n$.
Let $u$ be a solution of \eqref{Beam} such that
\begin{equation*}
\Vert u\Vert_{S(\mathbb{R})}<+\infty
\end{equation*}
then $u$ scatters.
\end{proposition}

\begin{proof}
We prove that $u$ scatters at $+\infty$. The proof for scattering at $-\infty$ is similar.
We claim that
\begin{equation}\label{AddedXu2}
\Vert \vert u\vert^{p-1}u\Vert_{L^\frac{2(n+2)}{n+4}(\mathbb{R}\times\mathbb{R}^n)}+\Vert \vert u\vert^{p-1}u\Vert_{L^\frac{2(n+4)}{n+8}(\mathbb{R}_+\times\mathbb{R}^n)}\lesssim_{E(u),\Vert u\Vert_S}1.
\end{equation}
First in case $1\le n\le 3$, we remark that
\begin{equation*}\label{AddedXuSmallDim}
\begin{split}
&\Vert \vert u\vert^{p-1}u\Vert_{L^\frac{2(n+2)}{n+4}(\mathbb{R}\times\mathbb{R}^n)}+\Vert \vert u\vert^{p-1}u\Vert_{L^\frac{2(n+4)}{n+8}(\mathbb{R}\times\mathbb{R}^n)}\\
&\lesssim \Vert u\Vert_S^{\theta p}\Vert u\Vert_{L^\infty(\mathbb{R}\times\mathbb{R}^n)}^{(1-\theta)p}+\Vert u\Vert_{S}^{\kappa p}\Vert u\Vert_{L^\infty(\mathbb{R}\times\mathbb{R}^n)}^{(1-\kappa)p}\\
&\lesssim_{E(u)} \Vert u\Vert_S^{\theta p}+\Vert u\Vert_{S}^{\kappa p}
\end{split}
\end{equation*}
for $\theta=(n+4)(p+1)/(2p(n+2))$ and $\kappa=(n+8)(p+1)/(2p(n+4))$. In the case $n=4$,
applying Strichartz estimates \eqref{Stric} and using conservation of Energy and H\"older's inequality, we get that
\begin{equation*}
\begin{split}
\Vert u\Vert_{L^{3}(\mathbb{R},L^{3p})}&\lesssim \sqrt{E(u)}+\Vert \vert u\vert^{p-1}u\Vert_{L^\frac{2(n+2)}{n+4}(\mathbb{R}\times\mathbb{R}^n)}+\Vert \vert u\vert^{p-1}u\Vert_{L^\frac{2(n+4)}{n+8}(\mathbb{R}\times\mathbb{R}^n)}\\
&\lesssim \sqrt{E(u)}+\Vert u\Vert_{S}^{p(1-\theta)}\Vert u\Vert_{L^\infty(\mathbb{R},L^{3p})}^{(p-1)\theta}\Vert u\Vert_{L^3(\mathbb{R},L^{3p})}^\theta\\
&+\Vert u\Vert_S^{p(1-\kappa)}\Vert u\Vert_{L^\infty(\mathbb{R},L^{3p})}^{(p-1)\kappa}\Vert u \Vert_{L^3(\mathbb{R},L^{3p})}^\kappa\\
&\lesssim_{E(u),\Vert u\Vert_S} 1+\Vert u\Vert_{L^3(\mathbb{R},L^{3p})}^{\theta}+\Vert u\Vert_{L^3(\mathbb{R},L^{3p})}^{\kappa}
\end{split}
\end{equation*}
for
$\theta=(p-2)/(2p-1)$, and $\kappa=3(p-3)/(4(2p-1))$.
This implies that
\begin{equation*}\label{ConditionOnVersion6}
\Vert u\Vert_{L^3(\mathbb{R},L^{3p})}<+\infty
\end{equation*}
 and boundedness of $u$ in the $S$-norm and H\"older's inequality give \eqref{AddedXu2} also in dimension $4$. Now, let $U^\sigma$ be such that $(U^\sigma(\sigma),U^\sigma_t(\sigma))=(0,0)$ and $U^\sigma_{tt}+\Delta^2U^\sigma+mU^\sigma= -\lambda\vert u\vert^{p-1}u$. Let also $w^\sigma$ be defined by
$$ (w^\sigma(t),w^\sigma_t(t))=\mathcal{W}(t+\sigma)(u(\sigma),u_t(\sigma)),$$
then for all $t$, $u(t)=w^\sigma(t)+U^\sigma(t)$, and applying Strichartz estimates on the interval $I=(\sigma,+\infty)$ gives with \eqref{AddedXu2} that
\begin{equation*}
\begin{split}
\Vert w^\sigma \Vert_{S(I)}&\le \Vert u\Vert_{S(I)}+\Vert U^\sigma\Vert_{S(I)}\\
&\lesssim \Vert u\Vert_{S(I)}+\Vert u\Vert_{L^\frac{2(n+4)p}{n+8}(I\times\mathbb{R}^n)}^p+\Vert u\Vert_{L^\frac{2(n+2)p}{n+4}(I\times\mathbb{R}^n)}^p\\
&\lesssim o(1)
\end{split}
\end{equation*}
as $\sigma\to +\infty$. Choosing $\sigma$ such that \eqref{LWPStatement} hold for $w^S$, we can apply Lemma \ref{CondForScat} and conclude that $u$ scatters.
\end{proof}

In our analysis, we need the following stability statement, the proof of which follows from repeated use of Strichartz estimates in the spirit of the proof of Proposition \ref{PropositionVersion6} above.

\begin{proposition}\label{Stability}
For any $E>0$, $E^\prime$ and $M>0$, there exists a constant $\delta_0>0$ such that for all $\delta\le \delta_0$,
for any $v\in \mathbb{E}_I$ an almost solution of \eqref{Beam} in the sense that
\begin{equation}\label{AlmostSol}
v_{tt}+\Delta^2v+mv-\lambda\vert v\vert^{p-1}v=e,
\end{equation}
for any
$(u_0,u_1)\in\mathcal{E}$, if $0\in I$ and
\begin{equation}\label{CondForStab}
\begin{split}
&E(u_0-v(0),u_1-v_t(0))\le E^\prime\\
&\Vert \pi_1\mathcal{W}(\cdot)(u_0-v(0),u_1-v_t(0))\Vert_{S(I)}\le \delta\\
&\Vert v\Vert_{S(I)}\le M\\
&\sup_{t\in I}E(v,v_t)\le E\\
&\Vert e\Vert_{N(I)}\le \delta
\end{split}
\end{equation}
then there exists a unique solution $u$ of \eqref{Beam} such that $(u(0),u_t(0))=(u_0,u_1)$. Besides, we have the following bound
\begin{equation}\label{StabCCl}
\begin{split}
&\Vert u-v\Vert_{S(I)}\le C\delta\\
&\Vert u\Vert_{L^\infty(I,H^2)}\le C\left(E+E^\prime+\delta\right)
\end{split}
\end{equation}
\end{proposition}

Note that one can replace the control on the $N$-norm of $e$ in \eqref{CondForStab} by the same control on the norm in \eqref{SecNotAlternativeNorms}.

\section{Asymptotic results and compactness}\label{CritSol}

As mentioned in the introduction, Theorem \ref{MainThm} is a consequence of the following theorem, where we establish that the $S$-norm is not only finite for strong solutions, but also bounded by a quantity depending only on the energy.

\begin{theorem}\label{Thm2}
Let $2\le n\le 4$, $\lambda>0$ and $p$ such that $1=8/n<p<\infty$. There exists a function $\Lambda$, given by \eqref{DefOfLambda} such that for all strong solutions $u$ of \eqref{Beam}, there holds that
\begin{equation}\label{UnifBoundOnSNormThm2}
\begin{split}
&\Vert u\Vert_{S(\mathbb{R})}\le \Lambda(E(u))
\end{split}
\end{equation}
where $E(u)$ is the energy defined in \eqref{SecNotEnengy}.
\end{theorem}

The proof of Theorem \ref{Thm2} will be completed in the next section.
As a remark, an analogue of Theorem \ref{Thm2} could be stated in the case $n \ge 5$. When $n \ge 5$ and
$1 + 8/n < p < 1 + 8/(n-4)$, we let the $S^\star$-norm be given by
\begin{equation*}\Vert u\Vert_{S^*(I)}=\Vert u\Vert_{L^\frac{2(n+2)p}{n+4}(I\times\mathbb{R}^n)}+\Vert u\Vert_{L^\frac{2(n+4)p}{n+8}(I\times\mathbb{R}^n)},
\end{equation*}
see Pausader [12]. Then, following the lines of the proof of Theorem \ref{Thm2}, we can prove that there exists a function $\Lambda$ such that for all strong solutions $u$ of (0.1), there holds that
$$\Vert u\Vert_{S^*(\mathbb{R})}\le  \Lambda(E(u)).$$

From now on, we focus on the proof of Theorem \ref{Thm2}. We always assume $1\le n\le 4$ and $1+8/n<p<\infty$. We also let $s=n(p-1)/(2(p+1))\in (0,2)$ be such that
\begin{equation}\label{DefOfsExponent}
\dot{H}^s\hookrightarrow L^{p+1}.
\end{equation}
In this section, we develop the study of sequences of solutions of \eqref{Beam} of bounded energy.

\begin{proposition}\label{SecAsympt1OscProp}
Let $(u_k)_k$ be a sequence of solutions of \eqref{Beam} such that \begin{equation}\label{CondForComp}
E(u_k)\to E_{max}\hskip.1cm\hbox{and}\hskip.1cm
\Vert u_k\Vert_{S}\to +\infty.
\end{equation}

Then one has
\begin{equation}\label{CCL1Osc}
\lim_{j\to+\infty}\limsup_{k\to +\infty}\sup_tE_0\left( u_k(t)-\sum_{\vert \log N\vert\le j} P_Nu_k(t)\right)=0.
\end{equation}
\end{proposition}

Following G\'erard \cite{Ger}, we say that the sequence $(u_k)_k$ is {\it uniformly $1$-oscillatory}. Informally speaking, this means that $u_k(t)$ remains localized around medium frequencies $\vert\xi\vert\sim 1$ for all times and uniformly in $k$.

\begin{proof}
Let $\omega$ be a solution of the linear equation \eqref{LinearBeam} such that $E_0(\omega)\le E_{max}$ and let $M$ be a dyadic number, then, using Strichartz estimates \eqref{Stric}, Bernstein's property \eqref{BernSobProp} and the isometric property of the linear flow $\mathcal{W}$, we get that, for $q=4/s=8(p+1)/(n(p-1))$,
\begin{equation}\label{HighFrequencyPartHasSmallNorm}
\begin{split}
\Vert P_{> M}\omega\Vert_{S}&\lesssim \Vert P_{>M}\omega\Vert_{L^\infty L^{p+1}}^{1-\frac{8}{n(p-1)}}\Vert P_{>M}\omega\Vert_{L^qL^{p+1}}^\frac{8}{n(p-1)}\\
&\lesssim \Vert P_{>M}\omega\Vert_{L^\infty \dot{H}^{s}}^{1-\frac{8}{n(p-1)}}\Vert \omega\Vert _{L^qL^{p+1}}^\frac{8}{n(p-1)}\\
&\lesssim E_0(\omega)^\frac{8}{n(p-1)}\left(\Vert P_{>M}\omega(0)\Vert_{\dot{H}^{s}}+\Vert P_{>M}\omega_t(0)\Vert_{\dot{H}^{s-2}}\right)^{1-\frac{8}{n(p-1)}}\\
&\lesssim E_{max}^\frac{8}{n(p-1)}\left( M^{2-s}E_{max}^\frac{1}{2}\right)^{1-\frac{8}{n(p-1)}}.\\
\end{split}
\end{equation}
Besides, still by \eqref{Stric} and \eqref{BernSobProp}, we also have that
\begin{equation}\label{SmallFrequencyPartHasSmallNorm}
\begin{split}
\Vert P_{<M}\omega\Vert_{S}&\lesssim \Vert P_{<M}\omega\Vert_{L^\infty L^\infty}^{1-\frac{2(n+4)}{n(p+1)}}\Vert P_{<M}\omega\Vert_{L^\frac{2(n+4)}{n}L^\frac{2(n+4)}{n}}^\frac{2(n+4)}{n(p+1)}\\
&\lesssim M^{\frac{n}{2}-\frac{n+4}{p+1}}E_0(\omega)^\frac{1}{2}\\
&\lesssim_{E_{max}} M^{\frac{n}{2}-\frac{n+4}{p+1}}.
\end{split}
\end{equation}
Consequently, we get that, for $M>1$ a dyadic number, there holds that
\begin{equation}\label{WeirdFrequencyPartHasSmallNorm}
\Vert P_{>M}\omega\Vert_S+\Vert P_{<M^{-1}}\omega\Vert_S\lesssim_{E_{max}} M^{-\eta}
\end{equation}
for some $\eta>0$.

Now, suppose that \eqref{CCL1Osc} is not true. Then there exists $\varepsilon>0$, and a sequence $j\to +\infty$ such that
\begin{equation}\label{Contrd}
\limsup_{k\to +\infty}\sup_tE_0\left(u_k(t)-\sum_{\vert \ln N\vert\le j}P_Nu_k(t)\right)>3\varepsilon.
\end{equation}
For any $j$ and any $k$, let $s=s(j,k)$ be such that
\begin{equation*}
E_0\left(u_k(s)-\sum_{\vert\ln N\vert\le j}P_Nu_k(s)\right)\ge \sup_t \left[E_0\left(u_k(t)-\sum_{\vert\ln N\vert\le j}P_Nu_k(t)\right)\right]-\varepsilon.
\end{equation*}
Then, for $j>8$, we choose two sequences $N_k\in (2^{-j},2^{-j}j)$ and $M_k\in (2^j/j,2^{j})$ of dyadic numbers such that
$$E_0(P^\prime_{N_k}u_k(s))\le 4E_{max}/\log j\hskip.5cm\hbox{and}\hskip.5cm E_0(P^\prime_{M_k}u_k(s))\le 4E_{max}/\log j,$$ where $P_N^\prime=P_{N/2}+P_N+P_{2N}$.
Let $\omega_{t_0,k}^0=\pi_1\mathcal{W}(\cdot)P_{N_k\le\cdot\le M_k}(u_k(t_0),\partial_tu_k(t_0))$, $\omega_{t_0,k}^{-}=\pi_1\mathcal{W}(\cdot)P_{< N_k}(u_k(t_0),\partial_tu_k(t_0))$ and $\omega_{t_0,k}^{+}=\pi_1\mathcal{W}(\cdot)P_{> M_k}(u_k(t_0),\partial_tu_k(t_0))$. 
We now claim that
\begin{equation}\label{SecAsymptSeparationEnergy4}
E(u_k)=E_0(\omega_{s,k}^{-})+E_0(\omega_{s,k}^{+})+E(\omega_{s,k}^0(t_0),\partial_t\omega_{s,k}^0(t_0))+o(1),
\end{equation}
where $o(1)\to 0$ as $j\to +\infty$.
Indeed,
\begin{equation}\label{SecAsymptSeparationEnergy1}
\begin{split}
\left\vert E_0(u_k(s))-E_0(\omega_{s,k}^{-})-E_0(\omega_{s,k}^0)-E_0(\omega_{s,k}^{+})\right\vert &\lesssim  E_0(P_{N_k}^\prime u_k(s))+E_0(P_{M_k}^\prime u_k(s))\\
&\lesssim_{E_{max}}1/\log j,
\end{split}
\end{equation}
while, by \eqref{BernSobProp}, for any $t\in\mathbb{R}$, we have that
\begin{equation}\label{SecAsymptSeparationEnergy2}
\Vert \omega_{s,k}^-(t)\Vert_{L^{p+1}}\lesssim N_k^s\Vert \omega_{s,k}^-(t)\Vert_{L^2}\lesssim_{E_{max}} (j2^{-j})^s.
\end{equation}
and, by Sobolev's inequality \eqref{DefOfsExponent},
\begin{equation}\label{SecAsymptSeparationEnergy3}
\begin{split}
\Vert \omega_{s,k}^+(t)\Vert_{L^{p+1}}&\lesssim \Vert \vert\nabla\vert^s\omega_{s,k}^+(t)\Vert_{L^2}\lesssim M_{k}^{s-2}\Vert \omega_{s,k}^+(t)\Vert_{H^2}\\
&\lesssim_{E_{max}} (2^j/j)^{s-2}.
\end{split}
\end{equation}
Now, \eqref{SecAsymptSeparationEnergy1}-\eqref{SecAsymptSeparationEnergy3} prove \eqref{SecAsymptSeparationEnergy4}.
Then, for $j$ sufficiently large, there holds that for $k$ sufficiently large,
$$E(P_{N_k\le\cdot\le M_k}(u_k(t_0),\partial_tu_k(t_0)))\le E_{max}-\varepsilon,$$ hence, letting $U_k$ be the nonlinear solution of \eqref{Beam} with initial data
$$(U_k(0),U_t(0))=P_{N_k\le\cdot\le M_k}(u_k(t_0),\partial_tu_k(t_0)),$$ we have the following bound by definition of $\Lambda$ in \eqref{DefOfLambda} and $E_{max}$ in \eqref{DefOfEmax}:
$$\Vert U_k\Vert_{S}\le \Lambda(E_{max}-\varepsilon)<+\infty.$$
Hence, using \eqref{WeirdFrequencyPartHasSmallNorm},  we have a sequence of solutions of \eqref{Beam} $U_k$ such that
\begin{equation}\label{Cond1OK}
\begin{split}
\Vert U_k\Vert_S\le M=\Lambda(E_{max}-\varepsilon)\hskip.5cm\hbox{and}\hskip.5cm E(U_k)\le E_{max}
\end{split}
\end{equation}
and a sequence of initial data $(u_k(s),\partial u_k(s))$ such that
\begin{equation}\label{Cond2OK}
\begin{split}
&E(u_k(s),\partial_tu_k(s))\le E_{max}\\
&\Vert \pi_1\mathcal{W}(\cdot)(U_k(0)-u_k(s),\partial_t(U_k(0)-u_k(s)))\Vert_{S}=o_j(1).
\end{split}
\end{equation}
Hence, applying proposition \ref{Stability}, we get that, for $j$ sufficiently large, and for all $k$ sufficiently large, there holds that
$u_k$ is global and that
$\Vert u_k\Vert_S\lesssim 1$, which contradicts \eqref{CondForComp}. This finishes the proof.
\end{proof}

Now we need the following structure result about sequences of $1$-oscillatory linear solutions of \eqref{LinearBeam}, Lemma \ref{SecAsymptLinearStructureLemma} below. In order to state it, we introduce some new definition.
We define a {\it core} to be any sequence $(S_k,Y_k)_k$ in $\mathbb{R}\times\mathbb{R}^n$ such that $(S_k)_k$ and $(Y_k)_k$ either go to infinity in norm or converge. Note that from any sequence in $\mathbb{R}\times\mathbb{R}^n$, one can always extract a core. We also say that two cores $(T_k,X_k)_k$ and $(S_k,Y_k)_k$ are {\it orthogonal} if
$$\vert T_k-S_k\vert+\vert X_k-Y_k\vert\to +\infty$$
as $k\to +\infty$. Our next lemma is inspired by the work in Bahouri and G\'erard \cite{BahGer} and M\'etivier and Schochet \cite{MetSch}. It reads as follows.

\begin{lemma}\label{SecAsymptLinearStructureLemma}
Let $(w_k)_k$ be a sequence of linear solutions of \eqref{LinearBeam} satisfying that $E_0(w_k)\le E$ for all $k$ and that $(w_k)_k$ is $1$-oscillatory\footnote{Note that for a linear solution, this is equivalent to $$\lim_{j\to+\infty}\limsup_{k\to +\infty}E_0(\omega_k(0)-\sum_{\vert\ln N\vert\le j}P_N\omega_k(0))=0$$}. Then there exists a sequence of linear solutions $(V^\alpha)_\alpha$ of \eqref{LinearBeam}, a sequence of cores $(S_k^\alpha,Y_k^\alpha)_\alpha$ which are pairwise orthogonal, and a sequence of linear solutions $(R^A_k)_{A,k}$ of \eqref{LinearBeam} such that, for all $k$, $A$, $t$ and $x$, there holds that
\begin{equation}\label{SecAsymptLinearStructureLemmaRes1}
\begin{split}
&w_k(t,x)=\sum_{\alpha=0}^AV^\alpha(t+S^\alpha_k,x+Y^\alpha_k)+R^A_k(t,x)\\
&E_0(w_k)=\sum_{\alpha=0}^AE_0(V^\alpha)+E_0(R^A_k)+o_A(1)\\
&\lim_{A\to +\infty}\limsup_{k\to +\infty}\Vert R^A_k\Vert_S=0\\
\end{split}
\end{equation}
where $o_A(1)\to 0$ for fixed $A$ as $k\to +\infty$.
Besides, for every $t\in\mathbb{R}$, the norm $\Vert w_k(t)\Vert_{L^{p+1}}$ converges as $k\to +\infty$, and there holds that
\begin{equation}\label{SecAsymptLinearStructureLemmaRes2}
\lim_{k\to +\infty}\Vert w_k(t)\Vert_{L^{p+1}}^{p+1}=\sum_{\alpha=0}^{\infty}\lim_{k\to +\infty}\Vert V^\alpha(t+S^\alpha_k)\Vert_{L^{p+1}}^{p+1}.
\end{equation}

\end{lemma}

\begin{proof}
Given a sequence of functions $W_k$ defined on $\mathbb{R}\times\mathbb{R}^n$, let us define
$P((W_k)_k)$ to be the set of weak limits (in $\mathcal{E}$) of sequences of the form $(\tau_{y_k}W_k(s_k),\partial_t\tau_{y_k}W_k(s_k))$,
where $(s_k,y_k)_k$ is an arbitrary core, and let
$$\mathcal{N}((W_k)_k)=\sup\{E_0(\phi,\psi):(\phi,\psi)\in P((W_k)_k)\}.$$

\medskip

Suppose first that $\mathcal{N}((w_k)_k)>0$, and let $(\phi,\psi)\in P((w_k)_k)$ be such that
$E_0(\phi,\psi)>\frac{1}{2}\mathcal{N}((w_k)_k)$. Then, we know that there exists $(S^0_k,Y^0_k)_k$ such that $(\tau_{Y^0_k}w_k(S^0_k),\partial_t\tau_{Y^0_{k}}w_k(S^0_k))\rightharpoonup (\phi,\psi)$. Let $W^0$ be the linear solution of equation \eqref{LinearBeam} with initial conditions $(\phi,\psi)$, and we let $w^1_k(t,x)=w_k(t,x)-\tau_{-Y^0_k}W^0(t-S^0_k,x)$.
By weak convergence, we have that
\begin{equation}\label{WeakEnergyDecouples}
\begin{split}
E_0(w_k)&=E_0(\tau_{Y^0_k}W^0(S^0_k))+E_0(\tau_{Y^0_k}w_k^1(S^0_k))+o(1)\\
&=E_0(W^0)+E_0(w^1_k)+o(1)
\end{split}
\end{equation}
where $o(1)\to 0$ as $k\to +\infty$.

\medskip

Then, we replace $P((w_k)_k)$ by $P((w^1_k)_k)$, and start the procedure again to define $(w^j_k)_k$ for $j=2,3,\dots$. If for some $j$, we find that $\mathcal{N}((w^j_k)_k)=0$, then we let $W^J=0$ and $Y^J_k=0$ for all $J> j$, and choose $(S^J_k)$ such that the cores $(S^J_k,0)$ are pairwise orthogonal and orthogonal to all cores $(S^p_k,Y^p_k)$, $0\le p\le j$.
Hence, we get a sequence $W^\alpha$ of solutions of \eqref{LinearBeam} and a sequence of cores $(S^j_k,Y^j_k)_k$. We first claim that the cores $(S^j_k,Y^j_k)_k$ are pairwise orthogonal.

\medskip

Suppose it is not so. Let $(S^i_k,Y^i_k)_k$ and $(S^j_k,Y^j_k)_k$ be such that $i< j$, $(S^i_k,Y^i_k)_k$ and $(S^j_k,Y^j_k)_k$ are not orthogonal but for every $0\le p<q<j$, there holds that $(S^p_k,Y^p_k)_k$ and $(S^q_k,Y^q_k)_k$ are orthogonal. In particular, we note that $\mathcal{N}((u^j_k)_k)>0$. Then passing to a subsequence, we may assume that $(S^j_k-S^i_k,Y^j_k-Y^i_k)\to (S_\ast,Y_\ast)$.  By definition, we have that
\begin{equation*}
\begin{split}
w^{j+1}_k(t,x)&=w_k(t,x)-\sum_{\alpha=0}^j W^\alpha(t-S^\alpha_k,x+Y^\alpha_k)\\
&=w^i_k(t,x)-\sum_{\alpha=i}^j W^\alpha(t-S^\alpha_k,x+Y^\alpha_k)\\
&=w^j_k(t,x)-W^j(t-S^j_k,x+Y^j_k),
\end{split}
\end{equation*}
and that for $l=i,j$, $(\tau_{Y^l_k}w^l_k(S^l_k),\partial_t\tau_{Y^l_k}w^l_k(S^l_k))$ converges to $(W^l(0),\partial_tW^l(0))$ weakly in $\mathcal{E}$.
Consequently, one the one hand, we find that
\begin{equation}\label{AddedEqt1}
\begin{split}
\tau_{Y^i_k}\bar{w}^i_k(S^i_k)-\tau_{Y^i_k}\bar{w}^j_k(S^i_k)&=\tau_{Y^i_k}\bar{w}^i_k(S^i_k)-\tau_{(-Y^\ast+o(1))}\tau_{Y^j_k}\bar{w}^j_k(S^j_k-S_\ast+o(1))
\end{split}
\end{equation}
converges weakly to
$\bar{W}^i(0)-\tau_{-Y_\ast}\bar{W}^j(-S_\ast)$, but the lefthand side of \eqref{AddedEqt1} above is also equal to
\begin{equation*}
\bar{V}=\sum_{\alpha=i}^{j-1}\tau_{Y^i_k}\bar{W}^\alpha(S^i_k-S^\alpha_k)
\end{equation*}
and since the cores $(S^\alpha_k,Y^\alpha_k)_k$ are pairwise orthogonal for $i\le \alpha\le j-1$, $\bar{V}$ converges weakly to
$\bar{W}^i(0)$. As a result, we get that $\bar{W}^j(-S_\ast)=0$, and hence $\bar{W}^j(0)=(0,0)$, but this gives a contradiction with the fact that $\mathcal{N}((w^j_k)_k)\ne 0$.
This proves that all the cores are pairwise orthogonal.

\medskip

Iterating \eqref{WeakEnergyDecouples}, we get that for any $A$ and $k$,
\begin{equation}\label{WeakEnergyDecouples2}
E_0(w_k)=\sum_{\alpha=0}^{A-1}E_0(W^\alpha)+E_0((w_k^A))+o(1).
\end{equation}
This shows that the series of $E_0(W^\alpha)$ is convergent, hence we get that $E_0(W^\alpha)\to 0$ as $\alpha\to +\infty$, and since $E_0(W^\alpha)>1/2\mathcal{N}((w_k^{\alpha+1})_k)$, this proves that
\begin{equation}\label{SecAsymptNgoesTo0}
\mathcal{N}((w_k^{\alpha})_k)\to 0\hskip.1cm\hbox{as}\hskip.1cm \alpha\to +\infty.
\end{equation}
Now, we need the following claim.

\medskip

For $N$ a dyadic number, there holds that
\begin{equation}\label{SecAsymptFixedFrequencyVanish}
\limsup_{k\to +\infty}\Vert P_Nw_k^A\Vert_{S}\lesssim_{E_{max}} N^\frac{n(p+1)-2(n+4)}{2(p+1)}\mathcal{N}((w_k)_k)^\frac{n(p+1)-2(n+4)}{2n(p+1)}. 
\end{equation}
To prove this, observe that, by Strichartz estimates \eqref{Stric}, we get that
\begin{equation}\label{AddedEqt2}
\Vert P_Nw_k^A\Vert_{L^\frac{2(n+4)}{n}(\mathbb{R},L^\frac{2(n+4)}{n})}\lesssim \sqrt{E_{max}}.
\end{equation}
Independently, we have the following estimation,
\begin{equation}\label{AddedEqt3}
\limsup_{k\to +\infty}\Vert P_Nw_k^A\Vert_{L^\infty (\mathbb{R}, L^\infty)}\lesssim N^\frac{n}{2}\mathcal{N}((w_k)_k)^\frac{1}{2}.
\end{equation}
Indeed, let $s_k$, and $x_k$ be such that $\Vert P_Nw_k^A\Vert_{L^\infty(\mathbb{R},L^\infty)}=\vert P_Nw_k^A(s_k,x_k)\vert$. Then, letting $\psi$ be as in \eqref{DefLitPalOp}, we get that
\begin{equation*}
\begin{split}
\limsup_{k\to +\infty}\Vert P_Nw_k^A(s_k,x_k)\Vert_{L^\infty(L^\infty)}
&\lesssim N^n\vert\int_{\mathbb{R}^n}w_k^A(s_k,y-x_k)\check{\psi}(Ny)dy\vert\\
&\lesssim N^n\sup_{(t_k)_k,(y_k)_k}\vert\int_{\mathbb{R}^n}w_k^A(t_k,y-y_k)\check{\psi}(Ny)dy\vert\\
&\lesssim N^n\sup_{W\in P((w_k^A)_k)}\vert \int_{\mathbb{R}^n}W(y)\check{\psi}(Ny)dy\vert\\
&\lesssim N^\frac{n}{2}\sup_{W\in P((w_k^A)_k)}\Vert W\Vert_{L^2}\\
&\lesssim N^\frac{n}{2}\mathcal{N}((w_k^A)_k)^\frac{1}{2}
\end{split}
\end{equation*}
This proves \eqref{AddedEqt3}, and \eqref{SecAsymptFixedFrequencyVanish} then follows from H\"older's inequality and \eqref{AddedEqt2}, \eqref{AddedEqt3}.

\medskip

Now, let us finish the proof of \eqref{SecAsymptLinearStructureLemmaRes1}. Using \eqref{WeakEnergyDecouples2} we see that we need only prove that
\begin{equation*}
\lim_{A\to +\infty}\limsup_{k\to +\infty}\Vert w^A_k\Vert_S=0.
\end{equation*}
But, using \eqref{SecAsymptFixedFrequencyVanish} and Strichartz estimates \eqref{Stric}, we get that
\begin{equation}\label{SecAsymptEstimForwA}
\begin{split}
\Vert w_k^A\Vert_S & \le \Vert P_{<L} w_k^A\Vert_S +\sum_{L\le N\le M}\Vert P_Nw_k^A\Vert_S+\Vert P_{>M}w_k\Vert_S\\
&\lesssim \sqrt{E_0(P_{<L}w_k)}+\sqrt{E_0(P_{>M}w_k)}+\sum_{L\le N\le M} N^\frac{n\theta}{2}\mathcal{N}((w_k)_k)^\frac{\theta}{2}\\
&\lesssim  \sqrt{E_0(P_{<L}w_k)}+\sqrt{E_0(P_{>M}w_k)}+\log\frac{M}{L} M^\frac{n\theta}{2}\mathcal{N}((w_k)_k)^\frac{\theta}{2}
\end{split}
\end{equation}
for $\theta=(n(p+1)-2(n+4))/(n(p+1))$.
Now, for every $\epsilon>0$, using \eqref{CCL1Osc} for $(w_k)_k$, we can choose $L$ and $M$ such that the sum of the two first term in the last line of \eqref{SecAsymptEstimForwA}
is smaller than $\epsilon$. Then \eqref{SecAsymptNgoesTo0} ensures that for $A$ sufficiently large, the third term is smaller than $\epsilon$.

\medskip

Finally, let us prove \eqref{SecAsymptLinearStructureLemmaRes2}. For simplicity, we will assume that $t=0$.
we remark that for $\alpha\ge 0$, there holds that
\begin{equation}\label{SecAsymptDecouplingLpNorm1}
\Vert \tau_{Y_k}V^\alpha(S^\alpha_k)+w^{\alpha+1}_k(0)\Vert_{L^{p+1}}^{p+1}=\Vert V^\alpha_k(S^\alpha_k)\Vert_{L^{p+1}}^{p+1}+\Vert w^{\alpha+1}_k(0)\Vert_{L^{p+1}}^{p+1}+o(1)
\end{equation}
Indeed, either $\vert S^\alpha_k\vert\to +\infty$, in which case, $\Vert V(S^\alpha_k)\Vert_{L^{p+1}}\to 0$, and \eqref{SecAsymptDecouplingLpNorm1} is clear, or there exists $S^\alpha_\ast$ such that $S^\alpha_k\to S^\alpha_\ast$. In this case, there holds that
\begin{equation*}
\begin{split}
&\left\vert\Vert \tau_{Y_k}V^\alpha(S^\alpha_k)+w^{\alpha+1}_k(0)\Vert_{L^{p+1}}^{p+1}-\Vert V^\alpha(S^\alpha_k)\Vert_{L^{p+1}}^{p+1}-\Vert w^{\alpha+1}_k(0)\Vert_{L^{p+1}}^{p+1}\right\vert\\
&\le \int_{\mathbb{R}^n}\vert \tau_{Y^\alpha_k}V^\alpha(S^\alpha_k)\vert \vert w^{\alpha+1}_k(0)\vert^pdx+\int_{\mathbb{R}^n}\vert\tau_{Y^\alpha_k}V^\alpha(S^\alpha_k)\vert^p\vert w^{\alpha+1}_k(0)\vert dx\\
& \le \int_{\mathbb{R}^n}\vert V^\alpha(S^\alpha_\ast)\vert \vert \tau_{-Y^\alpha_k}w^{\alpha+1}_k(0)\vert^pdx+\int_{\mathbb{R}^n}\vert V^\alpha(S^\alpha_\ast)\vert^p\vert \tau_{-Y^\alpha_k}w^{\alpha+1}_k(0)\vert dx+o_\alpha(1)\\
&\le o_\alpha(1),
\end{split}
\end{equation*}
where we used the fact that $\tau_{-Y^\alpha_k}w^{\alpha+1}_k(0)\rightharpoonup 0$ in $H^2$.
Besides, a similar proof shows that, for $\alpha\ne \beta$,
\begin{equation}\label{SecAsymptDecouplingLpNorm2}
\begin{split}
&\Vert \tau_{Y^\alpha_k}V^\alpha(S^\alpha_k)+\tau_{Y^\beta_k}V^\beta(S^\beta_k)\Vert_{L^{p+1}}^{p+1}\\
&=\Vert V^\alpha(S^\alpha_k)\Vert_{L^{p+1}}^{p+1}+\Vert V^\beta(S^\beta_k)\Vert_{L^{p+1}}^{p+1}+o_A(1).
\end{split}
\end{equation}
Let now define $\mathcal{F}$ to be the set of $\alpha$'s such that $(S^\alpha_k)_k$ is convergent to a limit $S^\alpha_\ast$. Equations \eqref{SecAsymptDecouplingLpNorm1} and \eqref{SecAsymptDecouplingLpNorm2} prove that, for every $A$, there holds that
\begin{equation*}
\begin{split}
\Vert w_k(t)\Vert_{L^{p+1}}^{p+1}&=\sum_{\alpha=0,\alpha\in \mathcal{F}}^A\Vert V^\alpha(S^\alpha_\ast)\Vert_{L^{p+1}}^{p+1}+\Vert w^{A+1}_k(0)\Vert_{L^{p+1}}^{p+1}+o(1)\\
&=\sum_{\alpha=0}^A\lim_{k\to +\infty}\Vert V^\alpha(S^\alpha_k)\Vert_{L^{p+1}}^{p+1}+\Vert w^{A+1}_k(0)\Vert_{L^{p+1}}^{p+1}+o(1),
\end{split}
\end{equation*}
where $o(1)\to 0$ as $k\to +\infty$.
To prove \eqref{SecAsymptLinearStructureLemmaRes2}, we need only show that
$$\limsup_{k\to +\infty}\Vert w^A_k\Vert_{L^\infty(\mathbb{R},L^{p+1})}\to 0\hskip.1cm\hbox{as}\hskip.1cm A\to +\infty.$$
Using \eqref{SecAsymptFixedFrequencyVanish}, we prove this as follows.
Let $t\in\mathbb{R}$, and let $s$ be as in \eqref{DefOfsExponent}, we have that
\begin{equation}\label{SecAsymptAddedEqtNormWA}
\begin{split}
\Vert w_k^A(t)\Vert_{L^{p+1}} & \le \Vert P_{<L}w_k^A(t)\Vert_{L^{p+1}}+\sum_{L\le N\le M}\Vert P_Nw_k^A(t)\Vert_{L^{p+1}}+\Vert P_{>M}w_k^A(t)\Vert_{L^{p+1}}\\
&\lesssim L^s\Vert w_k^A(t)\Vert_{L^2}+\Vert\nabla\vert^sP_{>M}w_k^A(t)\Vert_{L^2}\\
&+\sum_{L\le N\le M}\Vert P_Nw_k^A(t)\Vert_{L^\infty(\mathbb{R},L^2)}^\frac{2}{p+1}\Vert P_N w_k^A(t)\Vert_{L^\infty(\mathbb{R},L^\infty)}^\frac{p-1}{p+1}\\
&\lesssim L^s\sqrt{E_{max}}+M^{s-2}\sqrt{E_{max}}+\log\frac{M}{L}E_{max}^\frac{1}{p+1}\mathcal{N}((w_k^A)_k)^\frac{p-1}{p+1}.
\end{split}
\end{equation}
It follows that, for every $\epsilon>0$, we can choose $L$ and $M$ such that the sum of the two first term in the last line of \eqref{SecAsymptAddedEqtNormWA} are smaller that $\epsilon$. Then using \eqref{SecAsymptNgoesTo0}, we get that for $A$ and $k$ sufficiently large, the $L^\infty L^{p+1}$-norm of $w_k^A$ is smaller than $2\epsilon$.
This ends the proof of Lemma \ref{SecAsymptLinearStructureLemma}.

\end{proof}

Now we can prove the main result of this section

\begin{proposition}\label{SecAsymptNonlinearStructureProp}
Suppose that $E_{max}<+\infty$.
Let $(u_k)_k$ be a sequence of solutions of the nonlinear equation \eqref{Beam} and $(t_k)_k$ be a sequence in $\mathbb{R}$ such that
\begin{equation}\label{SecAsymptNonlinearStructurePropCond}
\begin{split}
E(u_k)\to E_{max}\hskip.3cm\hbox{and}\hskip.4cm\Vert u_k\Vert_{S(-\infty,t_k]},\Vert u_k\Vert_{S[t_k,+\infty)}\to +\infty\hskip.1cm,\hskip.1cm\hbox{as}\hskip.1cm k\to +\infty
\end{split}
\end{equation}
then there exists $(u_0,u_1)$ and a sequence $(y_k)_k$ in $\mathbb{R}^n$such that
\begin{equation}\label{SecAsymptNonlinearStructurePropCcl}
(\tau_{y_k}u_k(t_k),\tau_{y_k}\partial_tu_k(t_k))\to (u_0,u_1)\hskip.3cm\hbox{in}\hskip.3cm\mathcal{E}.
\end{equation}
\end{proposition}

\begin{proof}
By time translation, we may assume that $t_k=0$ for all $k$.
We know from Proposition \ref{SecAsympt1OscProp} that $(u_k)_k$ is $1$-oscillatory. Let us apply Lemma \ref{SecAsymptLinearStructureLemma} to the sequence of linear solutions
$(v_k)_k$ with initial data $(v_k(0),\partial_tv_k(0))=(u_k(0),\partial_tu_k(0))$. Then, we get
a sequence $(V^\alpha)_\alpha$ of \eqref{LinearBeam}, a sequence of cores $(S_k^\alpha,Y_k^\alpha)_\alpha$ which are pairwise orthogonal, and a sequence of linear solutions $(R^A_k)_{A,k}$ of \eqref{LinearBeam} such that, \eqref{SecAsymptLinearStructureLemmaRes1} holds true. Now, for all $\alpha$ define $U^\alpha$ as the nonlinear solution of \eqref{Beam} such that
$$\Vert V^\alpha(S^\alpha_k)-U^\alpha(S^\alpha_k)\Vert_{H^2}+\Vert \partial_tV^\alpha(S^\alpha_k)-\partial_tU^\alpha_k(S^\alpha_k)\Vert_{L^2}\to 0$$
as $k\to +\infty$. In particular, there holds that
\begin{equation}\label{SecAsymptNonlinearStructurePropProofEqt1}
\begin{split}
E(U^\alpha)&=E_0(V^\alpha)+\frac{\lambda}{p+1}\lim_{k\to +\infty}\Vert V^\alpha(S^\alpha_k)\Vert_{L^{p+1}}^{p+1}\\
\end{split}
\end{equation}
and that $E(U^\alpha)=0$ if and only if $E_0(V^\alpha)=0$.
Using \eqref{SecAsymptLinearStructureLemmaRes1}, \eqref{SecAsymptLinearStructureLemmaRes2} and \eqref{SecAsymptNonlinearStructurePropProofEqt1}, we get that
for fixed $A$,
\begin{equation*}
\begin{split}
E(u_k(0))&=E_0(\bar{u}_k(0))+\frac{\lambda}{p+1}\Vert u_k(0)\Vert_{L^{p+1}}^{p+1}\\
&=\sum_{\alpha=0}^A E_0(V^\alpha)+E_0(R^A_k)+o_A(1)+\frac{\lambda}{p+1}\sum_{\alpha=0}^\infty \lim_{k\to +\infty}\Vert V^\alpha(S^\alpha_k)\Vert_{L^{p+1}}^{p+1}\\
&=\sum_{\alpha=0}^A E(U^\alpha)+E_0(R^A_k)+o_A(1)+\frac{\lambda}{p+1}\sum_{\beta=A+1}^\infty \lim_{k\to +\infty}\Vert V^\alpha(S^\alpha_k)\Vert_{L^{p+1}}^{p+1}.
\end{split}
\end{equation*}
Letting $k\to +\infty$, and then $A\to +\infty$, we then get that
\begin{equation}\label{SecAsymptDecouplingEnergyForNonlinear}
E_{max}=\sum_{\alpha=0}^\infty E(U^\alpha)+\lim_{A\to +\infty}\limsup_{k\to +\infty}E_0(R^A_k).
\end{equation}
Suppose first that for all $\alpha$ there holds that $U^\alpha=0$, then, using \eqref{SecAsymptLinearStructureLemmaRes1}, we get that
\begin{equation*}
\Vert v_k\Vert_S\to 0
\end{equation*}
as $k\to +\infty$. Let $k$ be sufficiently large so that $\Vert v_k\Vert_S\le\delta$, we then get using Lemma \ref{CondForScat} that for all $k$ sufficiently large
$\Vert u_k\Vert_S\lesssim \Vert v_k\Vert_S\to 0$ which contradicts \eqref{SecAsymptNonlinearStructurePropCond}. Consequently, there exists $\alpha=\alpha_0$ such that $E(U^{\alpha_0})>0$. Suppose now that $E(U^{\alpha_0})<E_{max}$. Then, there exists $\epsilon>0$ such that for all $\beta$, $E(U^\beta)<E_{max}-\epsilon$. In particular, $\Vert U^\beta\Vert_S^{p+1}\le \Lambda(E_{max}-\epsilon)<+\infty$.

\medskip

In this case, let
$Q^A_k(s)=\sum_{\alpha=0}^A\tau_{Y^\alpha_k}U^\alpha(s+S^\alpha_k)$. We prove now that $Q^A_k$ is a good approximation of $u_k$. In the sequel, we define $U^\alpha_k$ by $U^\alpha_k(s,x)=U^\alpha(s+S^\alpha_k,x-Y^\alpha_k)$.

\medskip

First, there holds that
\begin{equation}\label{Claim1}
\Vert Q^A_k\Vert_{S}^{p+1}=\sum_{\alpha=1}^A\Vert U^\alpha\Vert_{S}^{p+1}+o_A(1)
\end{equation}
Indeed, we have that
\begin{equation}\label{SecAsympto1}
\begin{split}
&\left\vert\Vert Q^A_k\Vert_{S}^{p+1}-\sum_{\alpha=0}^A\Vert U^\alpha\Vert^{p+1}_{S}\right\vert\\
&=\left\vert\int_{\mathbb{R}\times\mathbb{R}^n}\left(\vert \sum_{\alpha=0}^AU^\alpha(t-S^\alpha_k,x-Y^\alpha_k)\vert^{p+1}-\sum_{\alpha=1}^A\vert U^\alpha(t-S^\alpha_k,x-Y^\alpha_k)\vert^{p+1}\right)dtdx\right\vert\\
&\lesssim\sum_{0\le \alpha\ne\beta\le A}\int_{\mathbb{R}\times\mathbb{R}^n}\vert U^\alpha(t-S^\alpha_k,x-Y^\alpha_k)\vert\vert U^\beta(t-S^\beta_k,x-Y^\beta_k)\vert^pdtdx\\
&=o_A(1).
\end{split}
\end{equation}
Now, since $\Lambda$ is sublinear around $0$ and bounded on $[0,E_{max}-\epsilon]$, using \eqref{SecAsymptDecouplingEnergyForNonlinear} and \eqref{Claim1}, we get that
\begin{equation*}
\begin{split}
\limsup_{k\to +\infty}\Vert Q^A_k\Vert_{S}^{p+1}&=\sum_{\alpha=0}^{A}\Vert U^\alpha\Vert_{S}^{p+1}\le\sum_{\alpha=0}^A\Lambda(E(U^\alpha))\lesssim_{E_{max},\epsilon} \sum_{\alpha=0}^AE(U^\alpha)\\
&\lesssim_{E_{max},\epsilon} E_{max}
\end{split}
\end{equation*}
Besides, $Q^A_k$ is an approximate solution in the sense that $Q^A_k$ satisfies \eqref{AlmostSol} with error
$$e=e^A_k=\sum_{\alpha=0}^Af(U^\alpha_k)-f(Q^A_k).$$
We claim that
\begin{equation}\label{IngredientLemma2}
\Vert e_k^A\Vert_{N}=o_A(1)
\end{equation}
as $k\to +\infty$.
To prove \eqref{IngredientLemma2}, we first prove that
\begin{equation}\label{SecAsymptLp+1Normofe}
\Vert e^A_k\Vert_{L^\frac{p+1}{p}(L^\frac{p+1}{p})}=o_A(1),
\end{equation}
where $o_A(1)\to 0$ as $k\to +\infty$ when $A$ is fixed. Indeed, proceeding as in \eqref{SecAsympto1}, we show that
\begin{equation*}
\begin{split}
\Vert e_k^A\Vert_{L^\frac{p+1}{p}(\mathbb{R},L^\frac{p+1}{p})}&\lesssim \Vert\sum_{\alpha\ne \beta}\vert U^\alpha_k\vert\vert U^\beta_k\vert^{p-1}\Vert_{L^\frac{p+1}{p}(\mathbb{R},L^\frac{p+1}{p})}\\
&\lesssim \sum_{\alpha\ne \beta}\Vert U^\alpha_k\vert U^\beta_k\vert^{p-1}\Vert_{L^\frac{p+1}{p}(\mathbb{R},L^\frac{p+1}{p})}\\
&\lesssim \sum_{\alpha\ne \beta}\Vert U^\alpha_k\vert U^\beta_k\vert^{p}\Vert_{L^1(\mathbb{R},L^1)}\lesssim o(1).
\end{split}
\end{equation*}
Then, if $n\le 3$, we let $\gamma=\eta=\infty$, $q=r=2(n+4)/n$, $\tau=\sigma=2(n+2)/n$ and $\theta=(p+1)/(pq^\prime)$, $\kappa=(p+1)/(p\tau^\prime)$, if $n=4$, we let
\begin{equation*}
\begin{split}
&\beta=\frac{1-\frac{(p+1)(n+8)}{2p(n+4)}}{\frac{\gamma(n+4)}{n}-\frac{p+1}{p}}\hskip.1cm,\hskip.1cm \theta=\frac{\frac{n+8}{2n}-\frac{1}{\gamma}}{\frac{(n+4)p}{n(p+1)}-\frac{1}{\gamma}}\hskip.1cm,\hskip.1cm\frac{1}{q}=\frac{n}{2(n+4)}+\beta\hskip.1cm,\hskip.1cm \frac{1}{r}=\frac{n}{2(n+4)}-\frac{4\beta}{n}\\
&\delta=\frac{1-\frac{(p+1)(n+4)}{2p(n+2)}}{\frac{\eta(n+2)}{n}-\frac{p+1}{p}}\hskip.1cm,\hskip.1cm \kappa=\frac{\frac{n+4}{2n}-\frac{1}{\eta}}{\frac{(n+2)p}{n(p+1)}-\frac{1}{\eta}}\hskip.1cm,\hskip.1cm\frac{1}{\tau}=\frac{n}{2(n+2)}
+\delta
\hskip.1cm\hbox{and}\hskip.1cm \frac{1}{\sigma}=\frac{n}{2(n+2)}-\frac{2\delta}{n}
\end{split}
\end{equation*}
with $\gamma$, $\eta$ large real numbers such that $2\le q,r,\tau,\sigma\le \infty$ and $0<\theta,\kappa<1$. On the one hand, using \eqref{SecAsymptLp+1Normofe} and the conservation of energy, we have that
\begin{equation*}
\begin{split}
\Vert e^A_k\Vert_{L^{q^\prime}(L^{r^\prime})}
&\le \Vert e^A_k\Vert_{L^\frac{p+1}{p}(L^\frac{p+1}{p})}^\theta\Vert e^A_k\Vert_{L^\infty(L^{\gamma})}^{1-\theta}\\
&\le o_A(1)\left(\Vert \sum_{\alpha=1}^A \tau_{Y^\alpha_k}U^\alpha_k\Vert_{L^\infty(H^2)}^p+\sum_{\alpha=1}^A\Vert U^\alpha\Vert_{L^\infty (H^2)}^p\right)^{1-\theta}\\
&\lesssim_{E_{max}} o_A(1),
\end{split}
\end{equation*}
and on the other hand, for the same reasons, there holds that
\begin{equation*}
\begin{split}
\Vert e^A_k\Vert_{L^{\tau^\prime}(L^{\sigma^\prime})}&\le \Vert e^A_k\Vert_{L^\frac{p+1}{p}(L^\frac{p+1}{p})}^\kappa \Vert e^A_k\Vert_{L^\infty(L^\eta)}^{1-\kappa}\\
&\le o_A(1)\left(\Vert \sum_{\alpha=1}^A \tau_{Y^\alpha_k}U^\alpha_k\Vert_{L^\infty(H^2)}^p+\sum_{\alpha=1}^A\Vert U^\alpha\Vert_{L^\infty (H^2)}^p\right)^{1-\kappa}\\
&\lesssim_{E_{max}} o_A(1).
\end{split}
\end{equation*}
With these inequalities, we get that \eqref{IngredientLemma2} holds true.

\medskip

Besides, by definition, there holds that
\begin{equation*}
\begin{split}
\Vert \pi_1\mathcal{W}(\cdot)(Q_k^A(0)-u_k(0),\partial_tQ_k^A(0)-\partial_tu_k(0))\Vert_{S}=\Vert \sum_{\alpha=0}^AW^\alpha-w_k\Vert_{S}=\Vert R^A_k\Vert_{S}\\
\end{split}
\end{equation*}
hence,
\begin{equation}\label{BlaBlaBlaBla}
\limsup_{k\to +\infty}\Vert \pi_1\mathcal{W}(\cdot)(Q_k^A(0)-u_k(0),\partial_tQ_k^A(0)-\partial_tu_k(0))\Vert_{S}\to 0
\end{equation}
as $A\to +\infty$.
Now, using \eqref{Claim1}, \eqref{IngredientLemma2}, \eqref{BlaBlaBlaBla} we can apply Proposition \ref{Stability} to get that
$$\Vert u_k\Vert_S\lesssim \lim_{A\to \infty}\limsup_{k\to +\infty}\Vert Q^A_k\Vert_{S}\lesssim_{E_{max},\epsilon}1,$$
which contradicts \eqref{SecAsymptNonlinearStructurePropCond}.

\medskip

Consequently, we find that for all $t$, $x$
\begin{equation*}
v_k(t,x)=V(t+S_k,x+Y_k)+R_k(t,x)
\end{equation*}
where $E_0(R_k)\to 0$ as $k\to +\infty$. Now, suppose that $ S_k\to +\infty$. Then, using the Strichartz estimates \eqref{Stric}, we get that
\begin{equation*}
\Vert v_k\Vert_{S([0,+\infty))}=\Vert V\Vert_{S([S_k,+\infty))}\to 0\hskip.1cm\hbox{as}\hskip.1cm k\to +\infty.
\end{equation*}
Then, using Lemma \ref{CondForScat}, we find that $\Vert u_k\Vert_{S([0,+\infty))}\lesssim \Vert v_k\Vert_{S([0,+\infty))}\to 0$. Again, this contradicts \eqref{SecAsymptNonlinearStructurePropCond}. A similar statement holds if $S_k\to -\infty$. Finally, there exists $S_\ast$ such that $S_k\to S_\ast$. Then we get that
\begin{equation*}
(u_k(0,x),\partial_tu_k(0,x))=(V(S_\ast,x+Y_k),\partial_tV(S_\ast,x+Y_k))+r_k
\end{equation*}
where $(r_k,\partial_tr^k)\to (0,0)$ in $\mathcal{E}$. This ends the proof of Proposition \ref{SecAsymptNonlinearStructureProp}.
\end{proof}

As a consequence of Proposition \ref{SecAsymptNonlinearStructureProp}, we have the following compactness result.

\begin{corollary}\label{ImportantCor}
Suppose that $E_{max}<+\infty$. Then there exists $u$, a strong solution of \eqref{Beam} defined on $\mathbb{R}$, and a function $y(t)$ such that the set
\begin{equation}\label{SecAsymptCompactnessRes}
K=\{(\tau_{y(t)}u(t),\tau_{y(t)}u_t(t)):t\in\mathbb{R}\}
\end{equation}
is precompact in $\mathcal{E}$. Besides, one can assume that $y$ is $C^1$ and satisfies
\begin{equation}\label{YLipschitz}
\vert \dot{y}(t)\vert\lesssim_u 1
\end{equation}
uniformly in $t$.
\end{corollary}

\begin{proof}

By definition of $E_{max}$, there exists a sequence $(u_k)_k$ of strong solutions of \eqref{Beam} satisfying the hypothesis of Proposition \ref{SecAsymptNonlinearStructureProp} with $t_k=0$ for all $k$. Let $(u_0,u_1)$ be given by Proposition \ref{SecAsymptNonlinearStructureProp}, and let $u$ be the solution of \eqref{Beam} with initial data $(u_0,u_1)$. Then,  $E(u)=E_{max}$, and by local wellposedness, we get\begin{equation*}
\Vert u\Vert_{S(-\infty,0)}=\Vert u\Vert_{S(0,+\infty)}=+\infty.
\end{equation*}
Now we need only find $y(t)$ satisfying the right properties.

\medskip

For $(u,v)\in\mathcal{E}$ and $y\in\mathbb{R}^n$, we define
\begin{equation*}
\begin{split}
&E_0(u,v,y,R)=\int_{B(y,R)}\left(v(x)^2+(\Delta u(x))^2+mu(x)^2\right)dx,\\
&\lambda(u,v,R)=\sup_{y\in\mathbb{R}^n}E_0(u,v,y,R)\hskip.3cm\hbox{and}\\
&\rho(u,v,\delta)=\inf\{R:\lambda(u,v,R)>(1-\delta)E_0(u,v)\}.
\end{split}
\end{equation*}
We claim that for all fixed $\delta>0$, $\rho(u(t),u_t(t),\delta)$ remains bounded. Indeed, if this were not true, there would exist a sequence of times $(t_k)_k$ such that we have for all $k$ and all $y$, that $E_0(u(t_k),u_t(t_k),y,k)\le (1-\delta)E_0(u(t_k),u_t(t_k))$ . But the sequence $(u(t_k),u_t(t_k))_k$ satisfies the hypothesis of Proposition \ref{SecAsymptNonlinearStructureProp}, hence there exists a sequence $Y_k$ such that, passing to a subsequence,
$$(\tau_{Y_k}u(t_k),\tau_{Y_k}u_t(t_k))\to (w_0,w_1)$$
as $k\to +\infty$ in $\mathcal{E}$. Consequently,
$$E_0(w_0,w_1,y,R)=\lim_k E_0(\tau_{Y_k}u(t_k),\tau_{Y_k}u_t(t_k),y,R)\le (1-\delta)E_0(w_0,w_1)$$
for all $R$. Hence $E_0(w_0,w_1)\le (1-\delta)E_0(w_0,w_1)$. This gives $E_0(w_0,w_1)=0$ which contradicts the fact that $E(w_0,w_1)=\lim E(u(t_k),u_t(t_k))=E_{max}$. Consequently, there exists an increasing function $R$ such that $\rho(u(t),u_t(t),\delta)<R(\delta)$ for all $t$.

A similar proof shows that there exists $\kappa(\delta)>0$ such that for all $t$,
\begin{equation}\label{SecAsymptdefinitionOfKappa}
\lambda(u(t),u_t(t),R(\delta))>\kappa(\delta).
\end{equation}
We now let $\delta$ be small such that $\delta<1/24$, and $\sqrt{\delta}<\kappa/(4E_{max})$.
Let $y(t)$ be such that $\lambda(u(t),u_t(t), R(\delta))=E(u(t),u_t(t),-y(t), R(\delta))$.

\medskip

We claim that the set $K=\{\tau_{y(t)}u(t),\tau_{y(t)}u_t(t):t\in\mathbb{R}\}$ is precompact in $\mathcal{E}$.
Suppose it were not so. Then there would exist $\epsilon>0$ and a sequence of times $t_i$ such that
\begin{equation}\label{SecAsymptNotCompactSequence}
E_0(\tau_{y(t_i)}u(t_i,)-\tau_{y(t_j)}u(t_j), \tau_{y(t_i)}u_t(t_i)-\tau_{y(t_j)}u_t(t_j))>\epsilon
\end{equation}
for all $i\ne j$. Applying Proposition \ref{SecAsymptNonlinearStructureProp}, we get that there exists a sequence $(Y_k)_k$ and $(w_0,w_1)\in\mathcal{E}$ such that, up to a subsequence,
$$(U(t_i),U_t(t_i))=(\tau_{y(t_i)+Y_i}u(t_i),\tau_{y(t_i)+Y_i}u_t(t_i))$$
converges to $(w_0,w_1)$ in $\mathcal{E}$.
In particular, $(U(t_i),U_t(t_i))$ is a Cauchy sequence. Let $i_0$ be such that for all $j\ge i_0$, there holds that
\begin{equation}\label{SecAsymptLastPropI0}
E_0(U(t_j)-U(t_{i_0}),U_t(t_j)-U_t(t_{i_0}))<\kappa/4,
\end{equation}
where $\kappa=\kappa(\delta)$ is defined in \eqref{SecAsymptdefinitionOfKappa},
and suppose that there exists a subsequence such that $\vert Y_k-Y_{i_0}\vert\to +\infty$ as $k\to +\infty$. Then, for $\vert Y_k-Y_{i_0}\vert>2R(\delta)$, we get that
\begin{equation}
\begin{split}
&E_0(U(t_j)-U(t_{i_0}),U_t(t_j)-U_t(t_{i_)}))\\
&=E_0(\bar{U}(t_j))+E_0(\bar{U}(t_{i_0}))+2\langle \bar{U}(t_j),\bar{U}(t_{i_0}\rangle\\
&\ge 2\kappa+2\int_{B(Y_j,R(\delta))}\left(U_t(t_j)U_t(t_{i_0})+\Delta U(t_j)\Delta U(t_{i_0})+mU(t_j)U(t_{i_0})\right)dx\\
&+2\int_{\vert x-Y_j\vert\ge R(\delta))}\left(U_t(t_j)U_t(t_{i_0})+\Delta U(t_j)\Delta U(t_{i_0})+mU(t_j)U(t_{i_0})\right)dx\\
&\ge 2\kappa-2\sqrt{E_0(\bar{U}(t_j))}\sqrt{\delta E_0(\bar{U}(t_{i_0}))}-2\sqrt{E_0(\bar{U}(t_{i_0}))}\sqrt{\delta E_0(\bar{U}(t_j))}\\
&\ge 2\kappa-4\sqrt{\delta}E_{max}\ge \kappa
\end{split}
\end{equation}
but this contradicts \eqref{SecAsymptLastPropI0}. Consequently, the sequence $(Y_k)_k$ remains bounded. Hence, up to a subsequence, we can assume that $Y_k\to Y_\ast$. But then we get that
$$(\tau_{-Y_\ast}U(t_k),\tau_{-Y_\ast}U_t(t_k))=(u(t_k,\cdot-y(t_k)-(Y_k-Y_\ast)),u_t(t_k,\cdot-y(t_k)-(Y_k-Y_\ast))))$$
is a Cauchy sequence, which contradicts \eqref{SecAsymptNotCompactSequence}. 

It only remains to prove \eqref{YLipschitz}. By the precompactness of $K$, and the continuity of the flow, there exists $s_0>0$ such that for every solution $u$ of \eqref{Beam} with initial data $(w_0,w_1)\in K$, there holds that
\begin{equation*}
\begin{split}
&E_0(u(s),u_t(s),0,2R(\delta))\ge (1-\delta)E_0(u(0),u_t(0),0,R(\delta))\hskip.1cm,\hskip.1cm\hbox{and}\\
&E_0(u(s),u_t(s))\le (1-\delta)^{-1}E_0(u(0),u_t(0))
\end{split}
\end{equation*}
for every time $s$ such that $\vert s\vert\le s_0$.
In particular,
$$E_0(u(s),u_t(s),0,2R(\delta))\ge (1-\delta)^3E_0(u(s),u_t(s)).$$ This implies that $E(u(s),u_t(s),Y,R(\delta))<(1-\delta)E_0(u(s),u_t(s))$ when $\vert Y\vert>3R(\delta)$. Consequently, for all $t\in\mathbb{R}$, for all sufficiently small $s \le 2s_0$, there holds that $\vert y(t)-y(t+s)\vert\le 6R(\delta)$. Now, let $t_j=js_0$ for $i\in\mathbb{Z}$, and  let $\tilde{y}$ be a smooth function such that $\tilde{y}(t_j)=y(t_j)$, and $\vert\tilde{y}^\prime(t)\vert\le 8R(\delta)s_0^{-1}$. Then $\vert y(t)-\tilde{y}(t)\vert\lesssim_u 1$, hence $\{\tau_{\tilde{y}(t)}u(t),t\in\mathbb{R}\}$ also is compact in $\mathcal{E}$ and replacing $y$ by $\tilde{y}$, we get \eqref{YLipschitz}.
This finishes the proof.
\end{proof}

Of course $u(t)$ might vanish for some $t$ (in which case the momentum vanishes). However, we prove now that it has a nontrivial contribution to the energy in every finite interval of time.

\begin{corollary}\label{NontrivialContribution}
Let $u$ be a nonlinear strong solution of \eqref{Beam} such that the set $K$ defined in \eqref{SecAsymptCompactnessRes} is precompact in $\mathcal{E}$, and $E(u)\ne 0$.
For every $\tau>0$, there exists two positive numbers $\alpha(\tau,u)$ and $\beta(\tau,u)$ such that, for all times $t$, there holds that
\begin{equation}\label{SectAsymptNontrivialNonlinearNorm}
\alpha\le \int_t^{t+\tau}\int_{\mathbb{R}^n}\vert u(s,x)\vert^{p+1}dxds\le \beta
\end{equation}
In particular, there holds that $\int_0^t\int_{\mathbb{R}^n}\vert u\vert^{p+1}dxds\simeq_u t$. 
\end{corollary}

\begin{proof}
The bound from above follows from Sobolev's inequality and the conservation of the energy. Suppose the bound from below is not true. Then there exist $\tau>0$ and a sequence $t_k$ such that
\begin{equation}\label{SectAsymptNontrivialNonlinearNormEqt1}
\int_{t_k}^{t_k+\tau}\int_{\mathbb{R}^n} \vert u(t,x)\vert^{p+1}dxdt<\frac{1}{k}.
\end{equation}
Using the precompactness of $K$, we can extract a subsequence and assume that
$$(\tau_{y(t_k)}u(t_k),\tau_{y(t_k)}u_t(t_k))\to (U_0,U_1)\hskip.1cm\hbox{in}\hskip.1cm \mathcal{E}.$$
Let $U$ be the nonlinear strong solution of \eqref{Beam} with initial data $(U_0,U_1)$. Then, $E(U)=E(u)\ne 0$. By wellposedness and \eqref{SectAsymptNontrivialNonlinearNormEqt1}, we get
\begin{equation*}
\int_0^\tau\int_{\mathbb{R}^n}\vert U(t,x)\vert^{p+1}dtdx=0.
\end{equation*}
Consequently, we have $U(t)=0$ for all $t$ in $(0,\tau)$, hence $U_t(t)=0$ for all such $t$. Consequently, $E(U)=0$. This is a contradiction.
\end{proof}

In the appendix, we give a more complete description of $y(t)$ in terms of $\hbox{Mom}(u)$.

\section{Proof of Theorems \ref{MainThm} and \ref{Thm2}}\label{ProofofThm}

The previous section showed that if $E_{max}<+\infty$, then there must exist a solution $u\ne 0$ of energy $E_{max}$ which is similar to a traveling wave. In order to rule out such a special solution, we need to understand how to rule out the existence of exact traveling waves for the defocusing equation. In the special case of a standing wave (i.e. traveling wave of speed $c=0$), we can use the usual Pohozhaev identity. 
In the special case of a solution similar to a traveling wave, but of momentum vector $0$\footnote{Note that a traveling wave has zero momentum vector if and only if it has speed $c=0$, i.e. if and only if it is a standing wave.},
a Virial estimate corresponding to the Pohozaev identity, (see Proposition \ref{SecPMTVirialProp}) allows us to get a contradiction. To prove that when $n\ge 2$, there can be no traveling wave with nonzero speed in the defocusing case, we use a variant of the Pohozhaev identity in a direction orthogonal to the momentum. In particular, we need to have at least two dimensions. The corresponding estimate adapted to solutions which are not given by an elliptic equation but nevertheless are similar to traveling waves is given by Proposition \ref{CorGeneral} below. We remark that in dimension $n=1$, a traveling wave of speed $c$, $u(t)=Q(x-ct)$ satisfies the equation
$$\Delta^2Q+c^2\Delta Q+mQ+\lambda\vert Q\vert^{-1}Q=0$$
which only consists of full derivatives, and whose linear part is {\it not} coercive for $c$ large.

\begin{proposition}\label{CorGeneral}
Let $2\le n\le 4$, $\lambda>0$ and $1+8/n<p<+\infty$.
Then there holds that $E_{max}=+\infty$.
\end{proposition}

\begin{proof}
Without loss of generality, one can assume that the momentum vector is parallel to the first vector in the canonical basis of $\mathbb{R}^n$. Thus
\begin{equation}\label{VanishingOfMomentum}
\int_{\mathbb{R}^n}u_t(t,x)\partial_ju(t,x)dx=0
\end{equation}
for $j\ge 2$.
Let $a(x)=\phi(x/R)$ where $\phi$ is a nonnegative smooth radially symmetric function such that $\phi=1$ on the ball of radius $1$ and $\phi$ is supported in the ball of radius $2$ and let $z=x-y(t)$. We define the modified Virial/Morawetz action as
\begin{equation}\label{DefI2}
I_2(t)=\int_{\mathbb{R}^n}z_2a(z)\partial_2u(t,x)u_t(t,x)dx,
\end{equation}
where $z_2=x_2-y_2(t)$ denotes the second component of $z$.
Integration by parts and \eqref{Beam} show that
\begin{equation}\label{TimeDerivativeOfI2}
\begin{split}
\partial_tI_2&=\int_{\mathbb{R}^n}\partial_t\left(z_2a(z)\right)u_t(t,x)\partial_2u(t,x)dx+\int_{\mathbb{R}^n}z_2a(z)\partial_2\frac{(u_t(t,x))^2}{2}dx\\
&-\int_{\mathbb{R}^n}z_2a(z)\partial_2u\left(\Delta^2u+mu+\lambda\vert u\vert^{p-1}u\right)dx\\
&=\frac{1}{2}\int_{\mathbb{R}^n}\left(-u_t^2+(\Delta u)^2+mu^2+\frac{2\lambda}{p+1}\vert u\vert^{p+1}\right)dx-2\int_{\mathbb{R}^n}\vert\nabla\partial_2u\vert^2dx\\
&+\dot{z}_2\int_{\mathbb{R}^n}u_t\partial_2u dx+\int_{\vert z\vert\ge R}\mathcal{O}_2\left(\vert u_t\vert^2+\vert u\vert^2+\vert\nabla u\vert^2+\vert\partial_{ij}u\vert^2\right)dx\\
\end{split}
\end{equation}
where, by \eqref{VanishingOfMomentum}, the third integral in the second equality vanishes, $\dot{z}=\frac{d}{dt}z$ and
\begin{equation}\label{DefOfRemainingI2}
\begin{split}
&\mathcal{O}_2\left(\vert u_t\vert^2+\vert u\vert^2+\vert\nabla u\vert^2+\vert\partial_{ij}u\vert^2\right)\\
&=\frac{1}{2}\left(z_2\partial_2a-(1-a)\right)\left(-u_t^2+(\Delta u)^2+mu+\frac{2\lambda}{p+1}\vert u\vert^{p+1}\right)\\
&+(1-a)\left(2\vert\nabla\partial_2u\vert^2-\dot{z}_2u_t\partial_2u\right)+\Delta a(\partial_2u)^2-z_2\Delta a\Delta u\partial_2u)\\
&-(z_2\dot{y}\cdot\nabla a(z))u_t\partial_2u-2(z_2\nabla a)\cdot\nabla \partial_2u\Delta u
\end{split}
\end{equation}
is bounded by a constant multiple of
$(u_t^2+\vert\partial_{ij}u\vert^2+\vert\nabla u\vert^2+u^2)$ and is supported on the set $\vert z\vert\ge R$. 
Besides, we define the equirepartition of energy action
$$J_a(t)=\int_{\mathbb{R}^n}a(z)u(t,x)u_t(t,x)dx.$$
Then
\begin{equation}\label{derivativeOfIa}
\begin{split}
\partial_tJ_a&=\int_{\mathbb{R}^n}a\left( u_t^2-(\Delta u)^2-mu^2-\lambda\vert u\vert^{p+1}\right)dx\\
&-\int_{\mathbb{R}^n}\left(\Delta a u\Delta u+2\nabla a\nabla u\Delta u\right)dx-\int_{\mathbb{R}^n}\left(\dot{y}\cdot\nabla a\right)u u_t\\
&=\int_{\mathbb{R}^n}\left(u_t^2-(\Delta u)^2-mu^2-\lambda\vert u\vert^{p+1}\right)dx\\
&-\int_{\mathbb{R}^n}\mathcal{O}_J(u)dx,
\end{split}
\end{equation}
where
\begin{equation}\label{DefOfRemainingJ}
\begin{split}
\mathcal{O}_J(u)=&(1-a)\left( u_t^2-(\Delta u)^2-mu^2-\lambda\vert u\vert^{p+1}\right)+\left(\Delta a u\Delta u+2\nabla a\nabla u\Delta u\right)\\
&+\left(\dot{y}\cdot\nabla a\right)u u_t
\end{split}
\end{equation}
has the same properties as $\mathcal{O}_2(u)$ in \eqref{DefOfRemainingI2}. Considering
$A_2=I_2+\frac{1}{2}J_a$, and combining \eqref{TimeDerivativeOfI2} and \eqref{derivativeOfIa}, we get that
\begin{equation}\label{BoundOnA}
\vert A_2(t)\vert\lesssim RE(u)
\end{equation}
for all time $t$
and
\begin{equation}\label{DerivativeOfA2}
\begin{split}
\partial_tA_2&=-2\int_{\mathbb{R}^n}\left(\vert\nabla\partial_2u\vert^2+\frac{\lambda (p-1)}{4(p+1)}\vert u\vert^{p+1}\right)dx+\int_{\vert z\vert\ge R}\left(\mathcal{O}_2(u)+\mathcal{O}_J(u)\right)dx.
\end{split}
\end{equation}
If $E_{max}<+\infty$, we integrate this identity from $0$ to $T>0$.
Using Corollary \eqref{NontrivialContribution}, we get that there exists $\alpha=\alpha(1,u)>0$ such that 
$$\int_0^T\int_{\mathbb{R}^n}\vert u(s,x)\vert^{p+1}dxds\ge \alpha T$$
for all $T>1$.
By compactness of $K$ proven in Corollary \ref{ImportantCor},
$$\frac{2(p+1)}{\lambda (p-1)}\int_{\vert z\vert\ge R}\left\vert\mathcal{O}_2(u)+\mathcal{O}_J(u)\right\vert dx\le \frac{\alpha}{2}$$
for $R$ sufficiently large.
Thus $-A_2(T)\gtrsim_uT$ for large $T$. This contradicts \eqref{BoundOnA}.

\end{proof}

As already mentioned, Corollary \ref{CorGeneral} finishes the proof of Theorems \ref{MainThm} and \ref{Thm2}.

\section{Appendix}

The goal of this appendix is to clarify the relation between the translation function $y(t)$ defined in \eqref{SecAsymptCompactnessRes} and the Momentum. It shows that $y(t)$ is roughly comparable to $t\hbox{Mom}(u)$.

\begin{proposition}\label{SecPMTVirialProp}
Let $u$ be a nonlinear strong solution of \eqref{Beam} such that the set $K$ as in \eqref{SecAsymptCompactnessRes} is precompact in $\mathcal{E}$, $\vert\dot{y}(t)\vert\lesssim 1$ and $E(u)\ne 0$. Then for every $\epsilon>0$, there exists a constant $C_\epsilon>0$ such that
\begin{equation}\label{SecPMTy(t)}
\left\vert y(t)\cdot Mom(u)+2\int_0^t\left(\Delta u(s)^2+\frac{\lambda n(p-1)}{4(p+1)}\vert u(s)\vert^{p+1}\right)ds\right\vert\le C_\epsilon+\epsilon t
\end{equation}
In particular, for $t>1$, $y(t)\cdot Mom(u)\simeq_u t$.
\end{proposition}

\begin{proof}
For $a$ and $z$ as in \eqref{DefI2}, let
\begin{equation*}\label{DefIa}
I_a(t)=\int_{\mathbb{R}^n}a(z)\left(z\cdot\nabla u(t,x)\right)u_t(t,x)dx
\end{equation*}
be the usual Morawetz/virial action. Integration by parts gives that
\begin{equation}\label{DerivativeOfI}
\begin{split}
\partial_tI_a&=\frac{n}{2}\int_{\mathbb{R}^n}\left(mu^2+(\Delta u)^2+\frac{2\lambda}{p+1}\vert u\vert^{p+1}-u_t^2\right)dx-2\int_{\mathbb{R}^n}(\Delta u)^2dx\\
&-\dot{y}\cdot \hbox{Mom}(u)\\
&+\int_{\mathbb{R}^n}\mathcal{O}_I(u)dx,
\end{split}
\end{equation}
where
\begin{equation}\label{DefOfRemainingIa}
\begin{split}
&\mathcal{O}_I(u)\\
&=-(1-a)\left(\frac{n}{2}\left(mu^2+(\Delta u)^2+\frac{2\lambda}{p+1}\vert u\vert^{p+1}-u_t^2\right)-2(\Delta u)^2-u_t\dot{y}\cdot\nabla u\right)\\
&-\frac{1}{2}\left(z\cdot\nabla a\right)\left(mu^2+(\Delta u)^2+\frac{2\lambda}{p+1}\vert u\vert^{p+1}-u_t^2\right)\\
&-\left(\Delta a\Delta u (z\cdot\nabla u)+2\partial_iaz_j\partial_{ij}u\Delta u+2\nabla a\nabla u\Delta u\right)-\left(\dot{y}\cdot\nabla a\right)\left( z\cdot\nabla u\right)u_t\\
\end{split}
\end{equation}
is bounded by a constant multiple of
$(u_t^2+\vert\partial_{ij}u\vert^2+\vert\nabla u\vert^2+u^2)$ and is supported on the set $\vert z\vert\ge R$.

Now, let $A=I_a+\frac{n}{2}J_a$ be the total action, we find that
\begin{equation}\label{DerivativeOfAction}
\begin{split}
\partial_tA&=-2\int_{\mathbb{R}^n}\left((\Delta u)^2+\frac{\lambda n(p-1)}{4(p+1)}\vert u\vert^{p+1}\right)dx-\dot{y}\cdot\hbox{Mom}(u)\\
&+\int_{\vert z\vert\ge R}\left(\mathcal{O}_I(u)+\frac{n}{2}\mathcal{O}_J(u)\right)dx.
\end{split}
\end{equation}
Integrating with respect to $t$, we get that
\begin{equation}\label{SecIntAction}
\begin{split}
&\left\vert y(t)\cdot\hbox{Mom}(u)+2\int_0^t\int_{\mathbb{R}^n}\left(\Delta u)^2+\frac{\lambda n(p-1)}{4(p+1)}\vert u\vert^{p+1}\right)dx\right\vert\\
&=\left\vert A(t)-A(0)-\int_{\mathbb{R}^n}\left(\mathcal{O}_I(u)+\frac{n}{2}\mathcal{O}_{J}(u)\right)dx\right\vert.
\end{split}
\end{equation}
Now, by compactness of $K$ in $\mathcal{E}$, where $K$ is as in \eqref{SecAsymptCompactnessRes}, we get that for every $\epsilon>0$, there exists $R(\epsilon)>0$ such that for all $v\in K$
$$\int_{\vert x\vert\ge R(\epsilon)}\left(\vert\mathcal{O}_I(v)\vert+\frac{n}{2}\vert\mathcal{O}_J(u)\vert\right) dx\le \epsilon.$$
Choosing $R=R(\epsilon)$ in \eqref{SecIntAction}  gives \eqref{SecPMTy(t)}. To prove the last assertion we use Corollary \eqref{SectAsymptNontrivialNonlinearNorm} and choose $\epsilon<\alpha(1,u)$ in \eqref{SecPMTy(t)}.

\end{proof}

\begin{corollary}\label{CorRadial}
If $u$ is as in the conclusion of Corollary \ref{ImportantCor}, then either $u=0$ or $\hbox{Mom}(u)\ne 0$. In particular, in the latter case $u$ cannot be radially symmetrical.
\end{corollary}

\begin{proof}
Indeed, if $u$ is radial, then $\hbox{Mom}(u)=0$. This contradicts Proposition \ref{SecPMTVirialProp}.
\end{proof}

The previous proposition described the component of $y(t)$ in the direction of $\hbox{Mom}(u)$. In the next proposition, we prove that in the other components, $y(t)$ does not evolve significantly.

\begin{proposition}
Suppose that $\hbox{Mom}(u)\ne 0$. Under the same conditions as in Proposition \ref{SecPMTVirialProp}, we have that
for all $\epsilon>0$, there exists a constant $C_\epsilon>0$ such that for any vector $Z$ orthogonal to $\hbox{Mom}(u)$, there holds that
\begin{equation}\label{OrthogonalPropOfx}
\left\vert Z\cdot y(t)\right\vert\le C_{\epsilon}+t\epsilon.
\end{equation}
\end{proposition}

\begin{proof}
For $i\ne j$, and $(u,v)\in \mathcal{E}$, let
\begin{equation*}
\omega_{ij}(u,v)=v(t,x)\left(\partial_iu(x)e_j-\partial_ju(x)e_i\right),
\end{equation*}
where $e_i$ denotes the $i$-th vector of the canonical basis in $\mathbb{R}^n$.
For $u$ a soultion of \eqref{Beam}, the quantity
$$\Omega_{ij}(u)=\int_{\mathbb{R}^n}\omega_{ij}(u(t),u_t(t))dx$$
is a conserved quantity, which we simply write $\Omega_{ij}$.
We consider the rotational action as follows
\begin{equation*}
R_{ij}=\int_{\mathbb{R}^n}a(z)\left(z\cdot\omega_{ij}(u(t),u_t(t))\right)dx
\end{equation*}
where $a$ and $z$ are as in the proof of Proposition \ref{SecPMTVirialProp} above.
Then, there holds that
\begin{equation*}
\begin{split}
\partial_tR_{ij}&=-\dot{y}\cdot \Omega_{ij}\\
&+\int_{\mathbb{R}^n}(1-a)\omega_{ij}dx-\dot{y}\cdot\int_{\mathbb{R}^n}\nabla a\left(z\cdot\omega_{ij}\right)dx\\
&-\frac{1}{2}\int_{\mathbb{R}^n}\left(\partial_iaz_j-\partial_jaz_i\right)\left(u_t^2-(\Delta u)^2-mu^2-\frac{\lambda\vert u\vert^{p+1}}{p+1}\right)dx\\
&-\int_{\mathbb{R}^n}\left(\Delta a\Delta u z\cdot\left(\partial_iue_j-\partial_jue_i\right)-2\Delta u\nabla a\cdot(\partial_iu e_j-\partial_jue_i)\right)\\
&-2\int_{\mathbb{R}^n}\partial_ka\left(\partial_{ik}ue_j-\partial_{jk}ue_i\right)\cdot zdx\\
&=-\dot{y}\cdot \Omega_{ij}+\int_{R\le \vert z\vert\le 2R}\mathcal{O}_\omega(u)dx.
\end{split}
\end{equation*}
Integrating this with respect to time gives
\begin{equation*}
y(t)\cdot\Omega_{ij}=R_{ij}(t)-R_{ij}(0)+\int_0^t\int_{\vert z\vert\ge R}\mathcal{O}_\omega(u)dxdt,
\end{equation*}
and this gives 
\eqref{OrthogonalPropOfx} since $\mathcal{O}_\omega(u)$ has the same properties as $\mathcal{O}_2(u)$ in \eqref{DefOfRemainingI2} and the $\Omega_{ij}$ span the orthogonal plan to $\hbox{Mom}(u)$.
\end{proof}

\medskip\noindent{\small ACKNOWLEDGEMENT:} The author expresses his deep thanks to Emmanuel Hebey and Walter Strauss for their
constant support and for stimulating discussions during the preparation of this work. He also wishes to thank Guixiang Xu for pointing out a mistake in Proposition \ref{PropositionVersion6}.

\end{document}